\documentclass{birkjour}
\usepackage{cite}
 \newtheorem{thm}{Theorem}

 \newtheorem{prop}[thm]{Proposition}
 \theoremstyle{definition}
 \newtheorem{defn}[thm]{Definition}
 \theoremstyle{remark}
 \newtheorem{rem}[thm]{Remark}
 
 \newtheorem{conj}[thm]{Conjecture}

\newcommand{\Hom}{\operatorname{Hom}}
\newcommand{\im}{\mathrm{im}}
\newcommand{\End}{\operatorname{End}}
\newcommand{\eqdef}{\stackrel{\text{def.}}{=}}
\def\bE{\mathbf{E}}

\def\Z{\mathbb{Z}}
\def\C{\mathbb{C}}
\def\H{\operatorname{\mathbb{H}}}
\def\rk{\operatorname{rk}}
\def\fd{\mathfrak{d}}
\def\md{\boldsymbol{\fd}}
\def\F{\operatorname{F}}

\def\ext{\text{ext}}

\newcommand{\be}{\begin{equation*}}
\newcommand{\ee}{\end{equation*}}

\newcommand \pd {{\partial}}
\newcommand \bpd {{\overline{\partial}}}

\newcommand{\id}{\operatorname{id}}
\newcommand{\Tr}{\operatorname{Tr}}
\newcommand{\tr}{\operatorname{tr}}
\newcommand{\str}{\operatorname{str}}
\def\rH{\operatorname{H}}
\def\cA{\mathcal{A}}
\def\cP{\mathcal{P}}
\def\cT{\mathcal{T}}
\def\cH{\mathcal{H}}
\def\cO{\mathcal{O}}
\def\Ob{\operatorname{Ob}}
\def\HF{\mathrm{HF}}
\def\PV{\operatorname{PV}}
\def\HPV{\operatorname{HPV}}
\def\Cob{\mathrm{Cob}}
\def\ioda{\boldsymbol{\iota}}
\def\bbpd{\boldsymbol{\bpd}}
\def\O{\mathrm{O}}

\def\0{{\hat{0}}}
\def\1{{\hat{1}}}
\def\i{\mathbf{i}}

\def\Jac{\operatorname{Jac}}

\def\PF{\mathrm{PF}}
\def\HPF{\mathrm{HPF}}

\def\DF{\operatorname{DF}}
\def\HDF{\operatorname{HDF}}

\def\vect{\mathrm{vect}}
\def\bdelta{\boldsymbol{\delta}}

\begin{document}

%
%
%
%
%
%
%
%
%

\title[A differential model for B-type LG theories]
 {A differential model for B-type\\ Landau-Ginzburg theories}

\author[E. M. Babalic]{E. M. Babalic}

\address{
	Center for Geometry and Physics\\
	Institute for Basic Science\\
	7 Cheongam-ro, Nam-gu, Pohang\\
	Republic of Korea 37673\\
	{\em and} \\
	Horia Hulubei National Institute for Physics and Nuclear Engineering \\
    Str. Reactorului no.30, P.O.BOX MG-6 \\
    RO-077125, Bucharest-Magurele, Romania}
\email{mirela@ibs.re.kr, mbabalic@theory.nipne.ro}

\author[D. Doryn]{D. Doryn}

\address{
	Center for Geometry and Physics\\
	Institute for Basic Science\\
	7 Cheongam-ro, Nam-gu, Pohang\\
	Republic of Korea 37673}
\email{dmitry@ibs.re.kr}

\author[C. I. Lazaroiu]{C. I. Lazaroiu}

\address{
	Center for Geometry and Physics\\
	Institute for Basic Science\\
	7 Cheongam-ro, Nam-gu, Pohang\\
	Republic of Korea 37673\\
	{\em and} \\
	Horia Hulubei National Institute for Physics and Nuclear Engineering \\
    Str. Reactorului no.30, P.O.BOX MG-6 \\
    RO-077125, Bucharest-Magurele, Romania}
\email{calin@ibs.re.kr}

\author[M. Tavakol]{M. Tavakol}

\address{
	Center for Geometry and Physics\\
	Institute for Basic Science\\
	7 Cheongam-ro,  Nam-gu, Pohang\\
	Republic of Korea 37673}

\email{mehdi@ibs.re.kr}




\begin{abstract}
We describe a mathematically rigorous differential model for B-type open-closed 
topological Landau-Ginzburg theories defined by a pair $(X,W)$, where $X$ is a non-compact
 K\"ahlerian manifold with holomorphically trivial canonical line bundle and $W$ is a
  complex-valued holomorphic function defined on $X$ and whose critical locus is 
  compact but need not consist of isolated points. 
  We also show how this construction specializes to the case when $X$ is Stein and
   $W$ has finite critical set, in which case one recovers a simpler mathematical model.
\end{abstract}

\maketitle

\section{Axiomatics of two-dimensional oriented open-closed TFTs}
  
Classical oriented open-closed topological LG (Landau-Ginzburg)
theories of type B are classical field theories defined on compact
oriented Riemann surfaces with corners and parameterized by pairs
$(X,W)$, where $X$ is a non-compact K\"ahlerian manifold and
$W:X\rightarrow \C$ is a non-constant holomorphic function defined on
$X$ and called the superpotential.  Previous work in the Mathematics literature 
assumed algebraicity of $X$ and $W$,
being mostly limited to very simple examples such as $X=\C^d$ and
generally assumed that the critical points of $W$ are isolated, in
which case topological D-branes can be described by matrix
factorizations. We do not impose such restrictions since there is no
Physics reason to do so. This leads to a much more general
description.\footnote{We use the results, notations
  and conventions of \cite{nlg1,nlg2}.  In our terminology
  ``off-shell'' refers to an object defined at cochain level while
  ``on-shell" refers to an object defined at cohomology level. The
  Physics constructions and arguments behind this work can be found in
  \cite{LG1,LG2}.}

A \textit{non-anomalous quantum oriented 2-dimensional open-closed
  topological field theory (TFT)} can be defined axiomatically
\cite{tft} as a symmetric monoidal functor from a certain symmetric
monoidal category $\Cob_2^\ext$ of labeled 2-dimensional oriented
cobordisms with corners to the symmetric monoidal category
$\vect^s_\C$ of finite-dimensional supervector spaces defined over~$\C$. The objects of the category $\Cob_2^\ext$ are finite disjoint
unions of oriented circles and oriented segments while the morphisms
are oriented cobordisms with corners between such, carrying
appropriate labels on boundary components (labels which can be identified
with the topological D-branes).  By definition, the {\em closed
  sector} of such a theory is obtained by restricting the monoidal
functor to the subcategory of $\Cob_2^\ext$ whose objects are disjoint
unions of circles and whose morphisms are ordinary cobordisms (without
corners).  It was shown in \cite{tft} that such a functor can be
described equivalently by an algebraic structure which we shall call a
{\em TFT datum}.  We start by describing certain simpler algebraic structures,
which form part of any such datum:

\begin{defn}
A pre-TFT datum  is an ordered triple
$(\cH,\cT,e)$ consisting of:
\begin{enumerate}
\item A finite-dimensional unital and supercommutative superalgebra
  $\cH$ defined over $\C$ (called the {\em bulk algebra}), whose unit
  we denote by $1_\cH$
\item A Hom-finite $\Z_2$-graded $\C$-linear category $\cT$ (called
  the {\em category of topological D-branes}), whose composition of
  morphisms we denote by $\circ$ and whose units we denote by $1_a\in
  \End_\cT(a)\eqdef\Hom_\cT(a,a)$ for all objects $a\in\Ob\cT$
\item A family $e=(e_a)_{a\in \Ob\cT}$ consisting of even
  $\C$-linear  \textit{bulk-boundary
    maps} $e_a:\cH\rightarrow \Hom_\cT(a,a)$ defined for each
  object $a$ of $\cT$
\end{enumerate}
such that the following conditions are satisfied:
\begin{description}
\item For any object $a\in \Ob\cT$, the map $e_a$ is a unital morphism
  of $\C$- superalgebras from $\cH$ to the endomorphism algebra
  $(\End_\cT(a),\circ)$
\item For any two objects $a,b\in \Ob\cT$ and for any $\Z_2$-homogeneous
  elements $h\in\cH$ and $t\in \Hom_\cT(a,b)$, we have:~ $e_b(h)\circ t=(-1)^{\deg h\,\deg t} t \circ e_a(h)$\,.
\end{description}
\end{defn}

\begin{defn}
A {\em Calabi-Yau supercategory of parity $\mu$} is a pair
$(\cT,\tr)$, where:
\begin{enumerate}
\itemsep 0.0em
\item $\cT$ is a $\Z_2$-graded and $\C$-linear Hom-finite category
\item $\tr=(\tr_a)_{a\in \Ob \cT}$ is a family of $\C$-linear maps
  $\tr_a:\End_\cT(a)\rightarrow \C$ of $\Z_2$-degree $\mu$
\end{enumerate} 
such that the following conditions are satisfied:
\begin{description}
\itemsep 0.0em
\item For any two objects $a,b\in \Ob\cT$, the $\C$-bilinear pairing
  $\langle \cdot ,\! \cdot \rangle_{a,b}\!:\! \Hom_\cT(a,b)\times
 \Hom_\cT(b,a)\rightarrow \C$ defined through:
\be
\langle t_1,t_2\rangle_{a,b}\eqdef \tr_b (t_1\circ t_2)~,~~\forall t_1\in \Hom_\cT(a,b)~,~\forall t_2\in \Hom_\cT(b,a)
\ee
is non-degenerate
\item For any $a,b\in \Ob\cT$ and any $\Z_2$-homogeneous
  elements $t_1\in \Hom_\cT(a,b)$ and $t_2\in \Hom_\cT(b,a)$, we
  have:
\be
\langle t_1,t_2\rangle_{a,b}=(-1)^{\deg t_1\,\deg t_2}\langle t_2,t_1\rangle_{b,a}~~.
\ee
\end{description}
If only the second condition above is satisfied, we say that $(\cT,\tr)$ is
a {\em pre-Calabi-Yau supercategory} of parity $\mu$.
\end{defn}

\begin{defn}
A {\em TFT datum} of parity $\mu$ is a system
$(\cH,\cT,e,\Tr,\tr)$, where:
\begin{enumerate}
\item $(\cH,\cT,e)$ is a pre-TFT datum
\item $\Tr:\cH\rightarrow \C$ is an even $\C$-linear map (called
  the {\em bulk trace})
\item $\tr=(\tr_a)_{a\in \Ob\cT}$ is a family of $\C$-linear maps
  $\tr_a:\End_\cT(a)\rightarrow \C$ of $\Z_2$-degree $\mu$ (called the {\em
  boundary traces})
\end{enumerate}
such that the following conditions are satisfied:
\begin{itemize}
\item $(\cH,\Tr)$ is a supercommutative Frobenius superalgebra. This
  means that the pairing induced by $\Tr$ on $\cH$ is
  non-degenerate
\item $(\cT,\tr)$ is a Calabi-Yau supercategory of parity $\mu$
\item The following condition (known as the {\em topological Cardy
  constraint}) is satisfied for all $a, b\in \Ob\cT$:
\be 
\Tr(f_a(t_1)f_b(t_2))=\str (\Phi_{ab}(t_1,t_2))~~,
~~\forall t_1\in \End_\cT(a)~,~\forall t_2\in \End_\cT(b)~.
\ee
Here, $\str$ denotes the supertrace on the finite-dimensional
$\Z_2$-graded vector space $\End_\C(\Hom_\cT(a,b))$ and:
\begin{description}
\item The $\C$-linear \textit{boundary-bulk} map
  $f_a:\End_\cT(a)\rightarrow \cH$ of $\Z_2$-degree $\mu$ is defined
  as the adjoint of the bulk-boundary map $e_a:\cH\rightarrow
 \End_\cT(a)$ with respect to the non-degenerate traces $\Tr$ and
  $\tr_a$:
\be
\Tr(h f_a(t))=\tr_a(e_a(h)\circ t)~~,~~\forall h\in \cH~,~\forall t \in \End_\cT(a)~~
\ee 
\item For any $a,b\in \Ob\cT$ and any $t_1\in\End_\cT(a)$ and
  $t_2\in \End_\cT(b)$, the $\C$-linear map
  $\Phi_{ab}(t_1,t_2):\Hom_\cT(a,b)\rightarrow \Hom_\cT(a,b)$ is
  defined through:
\be
\Phi_{ab}(t_1,t_2)(t)\eqdef t_2\circ t\circ t_1~~,~~\forall t\in \Hom_\cT(a,b)~~.
\ee
\end{description}
\end{itemize}
\end{defn}

\section{B-type open-closed Landau-Ginzburg theories}

\begin{defn}
A {\em Landau-Ginzburg pair} of dimension $d$ is a pair $(X,W)$ where:
\begin{enumerate}
\item $X$ is a non-compact K\"{a}hlerian manifold of complex dimension $d$
  which is \textit{Calabi-Yau} in the sense that the canonical line
  bundle $K_X\eqdef \wedge^d T^\ast X$ is holomorphically trivial.
\item $W:X\rightarrow \C$ is a {\em non-constant} complex-valued holomorphic
  function.
\end{enumerate}
The {\em signature} $\mu(X,W)$ of a Landau-Ginzburg pair $(X,W)$ is
defined as the mod 2 reduction\footnote{We denote by $\hat{k}\in \Z_2$ the mod 2 
reduction of any integer $k\in\Z$.} of the complex dimension $d$ of $X$.
\end{defn}
\noindent The {\em critical set of $W$} is defined as the set
of critical points of $W$: 
\be
Z_W\eqdef \{p\in X|(\pd W)(p)=0\}~.
\ee  

\subsection{The off-shell bulk algebra}

Let $(X,W)$ be a Landau-Ginzburg pair with $\dim_\C X=d$.  The
\textit{space of polyvector-valued forms} is defined through:
\be
\PV(X)=\bigoplus_{i=-d}^0 \bigoplus_{j=0}^d \PV^{i,j}(X)=
\bigoplus_{i=-d}^0 \bigoplus_{j=0}^d \cA^j(X, \wedge^{|i|} TX)~~,
\ee
where $\cA^j(X, \wedge^{|i|} TX)\equiv \Omega^{0,j}(X,\wedge^{|i|}
TX)$.  The \textit{twisted Dolbeault differential} induced by $W$ is
defined through $\delta_W\eqdef \bbpd+\ioda_W :\PV(X)\rightarrow
\PV(X)$, where $\bbpd$ is the  Dolbeault operator of
$\wedge TX$ (which satisfies $\bbpd(\PV^{i,j}(X))\subset
\PV^{i,j+1}(X)$), while $\ioda_W\eqdef-\i({\partial W})\lrcorner$ is
the contraction with the holomorphic 1-form $-\i \partial
W\in\Gamma(X,T^\ast X)$ (which satisfies $\ioda_W(\PV^{i,j}(X))\subset
\PV^{i+1,j}(X)$). 
Here $\i$ denotes the imaginary unit.
 Notice that $(\PV(X),\bbpd,\ioda_W)$ is a bicomplex
since $\bbpd^2=\ioda_W^2=\bbpd\ioda_W+\ioda_W\bbpd=0$.

\begin{defn}
The \textit{twisted Dolbeault algebra of polyvector-valued forms} of
the LG pair $(X,W)$ is the supercommutative $\Z$-graded $\O(X)$-linear
dg-algebra $(\PV(X),\delta_W)$, where $\PV(X)$ is endowed with the
total $\Z$-grading.
\end{defn}

\begin{defn}
The \textit{cohomological twisted Dolbeault algebra} of
$(X,W)$ is the supercommutative $\Z$-graded $\O(X)$-linear algebra
defined through:
\be
\HPV(X,W)= \rH(\PV(X),\delta_W)~~.
\ee
\end{defn}

\subsubsection{An analytic model for the off-shell bulk algebra}

\begin{defn}
The \textit{sheaf Koszul complex} of $W$ is the following complex of
locally-free sheaves of $\cO_X$-modules\footnote{We denote by $\cO_X$ 
the sheaf of holomorphic functions on $X$
and by $\O(X)=\Gamma(X,\cO_X)$ the ring of holomorphic functions on $X$. Here $\Gamma$ denotes
taking holomorphic sections, while $\Gamma_\infty$ denotes taking smooth sections.}:
\be
(Q_W):~~0 \to \wedge^d TX \stackrel{\ioda_W}{\rightarrow}  \wedge^{d-1} TX \stackrel{\ioda_W}{\rightarrow} \dots \stackrel{\ioda_W}{\rightarrow} \cO_X \to 0~~,
\ee
where $\cO_X$ sits in degree zero and we identify the exterior power
$\wedge^k T X$ with its locally-free sheaf of holomorphic sections.
\end{defn}

\begin{prop}
Let $\H(Q_W)$ denote the hypercohomology of the Koszul complex
$Q_W$. There exists a natural isomorphism of $\Z$-graded
$\O(X)$-modules:
\be
\HPV(X,W)\cong_{\O(X)} \H(Q_W)~~,
\ee
where $\HPV(X,W)$ is endowed with the total $\Z$-grading. Thus: 
\be
\rH^k(\PV(X),\delta_W)\cong_{\O(X)} \H^k(Q_W),~\forall k\in \{-d,\dots, d\}~~.
\ee
Moreover, we have
$\H^k(Q_W)=\bigoplus_{i+j=k}\bE_\infty^{i,j}$,
where $(\bE_r^{i,j},\mathbf{d}_r)_{r\geq 0}$ is a spectral sequence which
starts with:
\be
\bE_0^{i,j}:=\PV^{i,j}(X)= \cA^j(X,\wedge^{|i|} TX),~~\mathbf{d}_0=\bbpd~~,~ (i=-d,\dots, 0,~j=0,\dots, d)~.
\ee
\end{prop}

\subsection{The category of topological D-branes}

\begin{defn}
A \textit{holomorphic vector superbundle} on $X$ is a $\Z_2$-graded
holomorphic vector bundle defined on $X$, i.e. a complex holomorphic
vector bundle~$E$ endowed with a direct sum decomposition~$E=E^\0\oplus E^\1$, where $E^\0$ and $E^\1$ are holomorphic
sub-bundles of $E$.
\end{defn}

\begin{defn}
A \textit{holomorphic factorization} of $W$ is a pair
$a=(E,D)$, where $E=E^\0\oplus E^\1$ is a holomorphic vector
superbundle on $X$ and $D\in \Gamma(X,\End^\1\!(\!E))$ is a holomorphic
section of the bundle $\End^\1\!(E)\!=\! \Hom(\!E^\0\!,\!E^\1\!)\oplus
\Hom(E^\1\!,\!E^\0)\subset \End(E)$ which satisfies the condition $D^2=W
\id_E$.
\end{defn}

\subsection{The full TFT data}

\begin{defn}
The \textit{twisted Dolbeault category of holomorphic factorizations} of
$(X,W)$ is the $\Z_2$-graded $\O(X)$-linear dg-category $\DF(X,W)$
defined as follows:
\begin{description}
\item The objects of $\DF(X,W)$ are the holomorphic factorizations of
  $W$
\item Given two holomorphic factorizations $a_1=(E_1,D_1)$ and
  $a_2=(E_2,D_2)$, the Hom spaces:
\be
\Hom_{\DF(X,W)}(a_1,a_2)\eqdef \cA(X,\Hom(E_1,E_2)) 
\ee
are endowed with the total $\Z_2$-grading and with the twisted
differentials $\bdelta_{a_1,a_2}\eqdef \bbpd_{a_1,a_2}+\md_{a_1,a_2}$,
where $\bbpd_{a_1,a_2}$ is the  Dolbeault differential of $\Hom(E_1,E_2)$, 
while $\md_{a_1,a_2}$ is the {\em defect differential}:
\be
\md_{a_1,a_2}(\rho\otimes f)=(-1)^{\rk \rho} \rho \otimes (D_2\circ f)-(-1)^{\rk \rho +\sigma(f)} \rho\otimes (f\circ D_1)
\ee
\item  The composition of morphisms
$\circ\!:\!\cA(X,\Hom(E_2,\!E_3))\times \cA(X,\Hom(E_1,\!E_2))$
$\rightarrow \cA(X,\Hom(E_1,E_3))$
is determined uniquely by the condition:
\be
(\rho\otimes f)\circ (\eta\otimes g)=(-1)^{\sigma(f) \, \rk \eta}(\rho\wedge \eta) \otimes (f\circ g)
\ee
for all pure rank forms $\rho,\eta\in \cA(X)$ and all pure
$\Z_2$-degree elements $f\in \Gamma_\infty(X,\Hom(E_2,E_3))$ and $g\in
\Gamma_\infty(X, \Hom(E_1,E_2))$, where $\sigma(f)$ is the degree of $f$.
\end{description}
\end{defn}
\noindent We have (omitting indices): 
$\bdelta^2=\bbpd^2=\md^2=\bbpd\circ\md+\md\circ \bbpd=0$\,.

\begin{defn}
The \textit{cohomological twisted Dolbeault category of holomorphic
  factorizations} of $(X,W)$ is the $\Z_2$-graded $\O(X)$-linear
category defined as the total cohomology category of $\DF(X,W)$:
\be
\HDF(X,W)\eqdef\rH(\DF(X,W))~.
\ee
\end{defn}

\begin{thm}
Suppose that the critical set $Z_W$ is compact. Then the cohomology
algebra $\HPV(X,W)$ of $(\PV(X),\delta_W)$ is finite-dimensional over
$\C$ while the total cohomology category $\HDF(X,W)$ of $\DF(X,W)$ is
Hom-finite over $\C$. Moreover, one can define\footnote{Explicit
  expressions for $\Tr,\tr,e$ can be found in \cite{nlg1}.} a
bulk trace $\Tr\!:\!\HPV(X,W)\!\rightarrow \!\C$, boundary traces
$\tr_{a_1,a_2}\!:\Hom_{\HPV(X,W)}(a_1,a_2)\!\rightarrow \!\C$ and
bulk-boundary maps $e_a\!:\!\HPV(X,W)\rightarrow \End_{\HDF(X,W)}(a)$ such
that the system:
\be
(\HPV(X,W),\HDF(X,W), \Tr, \tr, e)
\ee
obeys the defining properties of a TFT datum except for non-degeneracy
of the bulk and boundary traces and for the topological Cardy
constraint.
\end{thm}

\begin{conj}
Suppose that the critical set $Z_W$ is compact. Then the quintuplet
$(\HPV(X,W), \HDF(X,W),\Tr, \tr, e)$ is a TFT datum and hence defines
a quantum open-closed TFT.
\end{conj}

\section{B-type Landau-Ginzburg theories on Stein manifolds}

When $X$ is a Stein manifold, Cartan's theorem B states that the
higher sheaf cohomology $\rH^i(X,\mathcal{F})$ vanishes when $i>0$
for any coherent analytic sheaf $\mathcal{F}$.

\subsection{An analytic model for the on-shell bulk algebra}
\begin{thm}
Suppose that $X$ is Stein. Then the spectral sequence defined
previously collapses at the second page $\bE_2$ and $\HPV(X,W)$ is concentrated in
non-positive degrees. 
For all $k=-d, \dots, 0$, the $\O(X)$-module
$\HPV^k(X,W)$ is isomorphic with the cohomology at position $k$ of the
following sequence of finitely-generated projective $\O(X)$-modules:
\be
(\cP_W):~~0 \to \rH^0(X, \wedge^d TX) \stackrel{\ioda_W}{\rightarrow}  \dots \stackrel{\ioda_W}{\rightarrow}  
\rH^0(X, TX) \stackrel{\ioda_W}{\rightarrow}  \O(X) \to 0~,
\ee
where $\O(X)$ sits in position zero. 
\end{thm}
\begin{rem}
Assume that $X$ is Stein. Then the critical set $Z_W$ is compact iff
it is finite, which implies $\dim_\C Z_W=0$\,.
\end{rem}
Let $\mathcal{J}_W\!\eqdef\! \im\,(\ioda_W\!\! :\! TX\to \cO_X)$ be the
critical sheaf of $W$,~ $Jac_W\eqdef \cO_X/\mathcal{J}_W$ be the Jacobi sheaf
  of $W$ and $\Jac(X,W)\eqdef \Gamma(X,Jac_W)$ be the Jacobi
  algebra of $(X,W)$\,.
\begin{prop}
Suppose that $X$ is Stein and $\dim_\C Z_W\!=\!0$. Then
$\HPV^k(X)\!=\!0$ for $k\neq 0$ and there exists a natural isomorphism
of $\O(X)$-modules:
\be
\HPV^0(X,W)\simeq_{\O(X)} \Jac(X,W)~~.
\ee
\end{prop}

\subsection{An analytic model for the category of topological D-branes}

\begin{defn}
The \textit{holomorphic dg-category of holomorphic
  factorizations} of $W$ is the $\Z_2$-graded $\O(X)$-linear 
  dg-category $\F(X,W)$ defined as follows:
\begin{description}
\item The objects are the holomorphic factorizations of $W$
\item Given two holomorphic factorizations $a_1=(E_1,D_1)$,
  $a_2=(E_2,D_2)$, let:
\be
\Hom_{\F(X,W)}(a_1,a_2)=\Gamma(X,\Hom(E_1,E_2))
\ee
be the space of morphisms between such, 
endowed with the $\Z_2$-grading: 
\be
\Hom_{\F(X,W)}^\kappa(a_1,a_2)=\Gamma(X, \Hom^\kappa(E_1,E_2))\,,~\forall \kappa\in \Z_2
\ee
and with the differentials $\md_{a_1,a_2}$ determined uniquely by the
condition:
\be
\quad\quad\md_{a_1,a_2}(f)\!=\! D_2\!\circ\! f \!-\!(-1)^\kappa f\!\circ\! D_1\,,~~\forall f\!\in \Gamma(X, \Hom^\kappa(E_1,\!E_2))\,,~\forall \kappa\in \Z_2
\ee
\item The composition of morphisms is the obvious one.
\end{description}
\end{defn}

\begin{thm}
Suppose that $X$ is Stein. Then $\HDF(X,W)$ is equivalent with the
total cohomology category $\HF(X,W)\eqdef \rH(\F(X,W))$ of $\F(X,W)$\,.
\end{thm}

\begin{defn}
A \textit{projective analytic factorization} of $W$ is a pair $(P,D)$,
where $P$ is a finitely-generated projective $\O(X)$-supermodule and
$D\in \End_{\O(X)}^\1(P)$ is an odd endomorphism of $P$ such that
$D^2=W \id_P$\,.
\end{defn}

\begin{defn}
The \textit{dg-category $\PF(X,W)$ of projective analytic
  factorizations} of $W$ is the $\Z_2$-graded $\O(X)$-linear
dg-category defined as follows:
\begin{description}
\item The objects are the projective analytic factorizations of $W$
\item Given two projective analytic factorizations $(P_1,D_1)$ and
  $(P_2,D_2)$ of $W$, set $\Hom_{\PF(X,W)}((P_1,D_1),(P_2,D_2))\eqdef
  \Hom_{\O(X)}(P_1,P_2)$ endowed with the obvious $\Z_2$-grading and
  with the $\O(X)$-linear odd differential
  $\md:=\md_{(P_1,D_1),(P_2,D_2)}$ determined uniquely by the
  condition:
\be
\md(f)= D_2\circ f-(-1)^{\deg f}f\circ D_1~,
~~\forall f\in \Hom_{\O(X)}(P_1,P_2)~
 \ee
\item The composition of morphisms is the obvious one.
\end{description}
\end{defn}

\begin{defn}
The \textit{cohomological category $\HPF(X,W)$ of analytic
  projective factorizations} of $W$ is the total cohomology category
$\HPF(X,W)\eqdef \rH(\PF(X,W))$, which is a $\Z_2$-graded $\O(X)$-linear
category.
\end{defn}

\begin{thm}
Suppose that $X$ is Stein. Then  $\HDF(X,W)$ and $\HPF(X,W)$ are
equivalent as $\O(X)$-linear $\Z_2$-graded categories.
\end{thm}
\noindent When $X$ is Stein and $Z_W$ is compact, the category of topological
D-branes of the B-type Landau-Ginzburg theory can be identified with
$\HPF(X,W)$.

\subsection*{Acknowledgment}
This work was supported by the research grant IBS-R003-S1.

\end{document}